\documentclass[11pt]{article}

\makeatletter
\@addtoreset{equation}{section}
\makeatother
\def\w{\omega}

\begin{document}

\def \beq{\begin{equation}}
\def\eeq{\end{equation}}
\newtheorem{theo}{Theorem}
\newtheorem{defi}{Definition}
\newtheorem{prop}{Proposition}

\title{The hyperbolic volume of knots from quantum dilogarithm}

\author{R.M. Kashaev\thanks{On leave of absence from
St. Petersburg Branch of the Steklov Mathematical Institute,
Fontanka 27, St. Petersburg 191011, RUSSIA}\\ \\
Laboratoire de Physique Th\'eorique 
{\sc enslapp}\thanks{URA 14-36 du CNRS,
associ\'ee \`a l'E.N.S. de Lyon, au LAPP d'Annecy 
et \`a l'Universit\`e de Savoie}
\\
ENSLyon,
46 All\'ee d'Italie,\\
69007 Lyon, FRANCE\\
E-mail: {\sf rkashaev\@@enslapp.ens-lyon.fr}}
\date{January, 1996}
\maketitle

\abstract{The invariant of a link in three-sphere, associated with the cyclic 
quantum dilogarithm, depends on a natural number $N$. 
By the analysis of particular examples it is argued that for a hyperbolic knot 
(link) the absolute value of this invariant  grows exponentially at large $N$, 
the hyperbolic volume of the knot (link) complement being the growth rate.}
\vskip 2cm

\rightline{{\small E}N{\large S}{\Large L}{\large A}P{\small P}-L-577/96}

\newpage

\section{Introduction}
 
Many known knot and link invariants, including Alexander  \cite{A} and Jones 
\cite{J} polynomials, can be obtained from 
$R$-matrices, solutions to the Yang-Baxter equation (YBE) \cite{Y,B1}.
Remarkably, many $R$-matrices in turn appear in quantum 3-dimensional 
Chern-Simons (CS) theory as partition functions 
of a 3-manifold with boundary, and with properly chosen Wilson lines \cite{W1}. The corresponding knot invariant acquires an interpretation in terms of a 
meanvalue of a Wilson loop. Note, that the quantum CS theory is an example of 
the topological quantum field theory (TQFT) defined axiomatically in \cite{At}.

 Thurston in his theory of hyperbolic 3-manifolds \cite{Th} introduces the 
notion of a hyperbolic knot: a knot that has a complement that can be given a 
metric of negative constant curvature. The volume of the complement in this 
metric
appears to be a topological invariant, called the hyperbolic volume of a knot 
\cite{Th,AHW,Ad}. In principle, 
there should exist a quantum generalization of this invariant for the following reason. 

Consider Euclidean quantum 
2+1 gravity with a negative cosmological constant, which is the CS theory with 
a non-compact gauge group \cite{W2}, and calculate the partition function of a 
hyperbolic knot's complement. The result would be a topological invariant, and 
the classical limit in the leading order would reproduce the hyperbolic volume 
of the complement.
Unfortunately, quantum 2+1 gravity in its present formulation can not
yet be used as a computational tool for this invariant.

The hyperbolic volume of ideal tetrahedra in three-dimensional hyperbolic space
can be expressed in terms of Lobachevsky's function, which is the imaginary 
part of Euler's dilogarithm \cite{M}. Therefore, it is natural to expect, that 
the quantum dilogarithm of \cite{FK,F1,F2} can lead to a generalized (deformed) notion of the hyperbolic volume.

In \cite{K1,K2} a link invariant, depending on a positive integer parameter $N$, has been defined via 3-dimensional interpretation of the cyclic quantum 
dilogarithm \cite{FK,BR}, see also \cite{BB,KMS}. The construction 
can be considered as an example of the simplicial (combinatorial) version of 
the 3-dimensional TQFT \cite{TV}. In this paper  we argue that this invariant 
is in fact a quantum generalization of the hyperbolic volume invariant. Namely, let $\langle L\rangle$ be the value of the invariant on a 
hyperbolic knot or link $L$ in three-sphere. We study the ``classical'' limit 
$N\to \infty$ of this invariant on particular examples of $L$, and show that 
for hyperbolic knots its absolute value grows exponentially, 
with the growth rate being 
given by the hyperbolic volume of the knot complement:
\beq\label{main}
2\pi\log|\langle L\rangle|\sim N V(L), \quad N\to\infty,
\eeq
where $V(L)$ is the hyperbolic volume of $L$'s complement
in $S^3$. Formula (\ref{main}) is in agreement with the expected classical 
limit of Euclidean quantum 2+1 gravity with a negative cosmological constant 
\cite{W2}. It is thus possible that the simplicial TQFT, defined in terms of 
the cyclic quantum dilogarithm, can be associated with quantum 2+1 dimensional 
gravity.

\section{The quantum invariant for three hyperbolic knots}\label{sec1}

Let $\w$ be a primitive $N$-th root of unity. Throughout this paper we will 
work with the following choice for this root:
\beq
\w=\exp(2\pi i/N).
\eeq
The result of calculation of the link invariant from papers \cite{K1,K2} for 
three simplest hyperbolic knots reads:
\beq\label{4_1}
\langle 4_1\rangle=\sum_{k}|(\w)_k|^2,\quad (\mbox{"figure-eight" knot}),
\eeq
\beq\label{5_2}
\langle 5_2\rangle=\sum_{k\le l}\frac{(\w)_l^2}{(\w)_k^*}\w^{-k(l+1)},
\eeq
\beq\label{6_1}
\langle 6_1\rangle=\sum_{k+l\le m}\frac{|(\w)_m|^2}{(\w)_k(\w)_l^*}
\w^{(m-k-l)(m-k+1)},
\eeq
where the summation variables run over $\{0,\ldots,N-1\}$;
\beq
(\w)_k=\prod_{j=1}^k(1-\w^j),\quad k=0,\ldots,N-1;
\eeq
and the asterisk means the complex conjugation. Strictly speaking only $N$-th 
powers of these quantities are invariants. 

Note that the simplest case
$N=2$ is related to a particular value of the Alexander polynomial 
$\Delta_L(t)$:
\beq
|\langle L\rangle |=\Delta_L(-1),\quad N=2.
\eeq
In the next section we study another extreme case $N\to\infty$.

\section{The classical limit}

Here we calculate explicitly the leading asymptotics at $N\to\infty$ of the 
invariant for the hyperbolic knots from section \ref{sec1} and justify formula
(\ref{main}).

First, for a positive real $\gamma$ and complex
$p$ with 
\beq
|\mbox{Re}\ p|<\pi+\gamma,
\eeq
define two functions
\beq\label{mainfun}
f_\gamma(p)=S_\gamma(\gamma-\pi)/S_\gamma(p), \quad
\overline f_\gamma(p)=S_\gamma(-p)/S_\gamma(\pi-\gamma),
\eeq
where
\beq\label{fad}
S_\gamma(p)=\exp\frac{1}{4}\int_{-\infty}^{+\infty}\frac{e^{px}}
{\sinh(\pi x)\sinh(\gamma x)}\frac{dx}{x},
\eeq
the singularity of the integrand at $x=0$ being put below the contour of 
integration.
In \cite{F1} the function (\ref{fad}) is shown to be a particular solution 
to the functional equation:
\beq\label{funeq}
(1+e^{ip})S_\gamma(p+\gamma)=S_\gamma(p-\gamma).
\eeq
Via this functional equation the definition of the function
$S_\gamma(p)$ can be extended to the whole complex plane. 

For fixed $p$ the leading asymptotics of $S_\gamma(p)$ at $\gamma\to 0$  
is given by Euler's dilogarithm:
\beq\label{as}
S_\gamma(p)\sim\exp\frac{1}{2i\gamma}\mbox{Li}_2(-e^{ip}),\quad \gamma\to0,
\eeq
where
\beq
\mbox{Li}_2(z)=-\int_0^z\frac{\log(1-u)}{u}du.
\eeq
For numerical calculations the following formula will be useful (see, for 
example, \cite{K}) :
\beq
\mbox{Im}\ \mbox{Li}_2(re^{i\theta})=
\varphi\log(r)+\Lambda(\varphi)+\Lambda(\theta)-\Lambda(\varphi+\theta),
\eeq
where $0<r\le1$,
\beq
\varphi=\varphi(r,\theta)=\arctan\left(\frac{r\sin\theta}{1-r\cos\theta}\right),
\eeq
and 
\beq\label{lob}
\Lambda(\theta)=-\int_0^\theta\log|2\sin\phi|\ d\phi,
\eeq
is Lobachevsky's function.

From (\ref{mainfun}) and (\ref{funeq}) it is easy to see that $f_\gamma(p)$ and 
$\overline f_\gamma(p)$ are analytic continuations of the symbols $(\w)_k$ and
$(\w)_k^*$ in the sense that
\beq\label{anal}
(\w)_k=f_\gamma(-\pi+\gamma+2k\gamma),\quad
(\w)_k^*=\overline f_\gamma(-\pi+\gamma+2k\gamma),
\eeq
where
\beq\label{gamma}
\quad\gamma=\pi/N.
\eeq
Formulae (\ref{anal}) enable us to rewrite the summations in (\ref{4_1}),
(\ref{5_2}), (\ref{6_1}) as contour integrals, one has just to replace
symbols $(\w)_k$ and $(\w)_k^*$ by their analytic continuations,
and each summation, by a contour integral:
\beq
\sum_k\to\frac{i}{4\gamma}\oint dp\
\tan\left(\frac{\pi^2+\pi p}{2\gamma}\right) 
\eeq
with an appropriately chosen contour.

In what follows $\gamma$ will be assumed to be specified as in (\ref{gamma}).

\subsection{ Knot $4_1$ (figure-eight knot)}

The exact formula (\ref{4_1}) can be rewritten as a contour integral:
\beq\label{contour4_1}
\langle 4_1\rangle=\frac{i}{4\gamma}\oint_Cdp\  
\tan\left(\frac{\pi^2+\pi p}{2\gamma}\right)f_\gamma(p)\overline f_\gamma(p),
\eeq
Where contour $C$ encircles $N$ points:
\beq
-\pi+\gamma+2k\gamma\quad 0\le k<N,
\eeq
in the counterclockwise direction. In the large $N$ (small $\gamma$) limit 
integral (\ref{contour4_1}), by using (\ref{as}), asymptotically can be 
approximated by
\beq
\langle 4_1\rangle\sim\int dz\ 
\exp\frac{i}{2\gamma}[\mbox{Li}_2(z)-\mbox{Li}_2(1/z)].
\eeq
The saddle point approximation of the last integral gives the result:
\beq\label{as41}
\langle 4_1\rangle\sim
\exp\frac{V(4_1)}{2\gamma},\quad \gamma=\pi/N\to0,
\eeq
where
\beq\label{vol4_1}
V(4_1)=4\Lambda(\pi/6)=2.02988321\ldots
\eeq
Formula (\ref{vol4_1}) is in agreement with the known hyperbolic volume of the 
figure-eight knot complement \cite{AHW}.

\subsection{Knot $5_2$}

The sum (\ref{5_2}) at large $N$ is approximated by the double contour integral
\beq\label{contour5_2}
\langle 5_2\rangle\sim\int dzdu\ 
\exp\frac{i}{2\gamma}[2\mbox{Li}_2(z)+\mbox{Li}_2(1/u)+\alpha(z,u)
-\pi^2/2],
\eeq
where
\beq\label{alpha}
\alpha(z,u)=\log(z)\log(u).
\eeq
The stationary points are described by the algebraic equations:
\beq\label{stat1}
u+z=uz,\quad u=(1-z)^2.
\eeq
The maximal contribution to the integral comes from the solution $(z_0,u_0)$ 
to (\ref{stat1}) with the property:
\beq
\mbox{Im}\ z_0<0<\mbox{Im}\ u_0.
\eeq
Thus, we obtain for asymtotics of the absolute value of the invariant:
\beq
|\langle 5_2\rangle|\sim \exp\frac{V(5_2)}{2\gamma},\quad \gamma=\pi/N\to0,
\eeq
where 
\beq
V(5_2)=-\mbox{Im}\ (2\mbox{Li}_2(z_0)+\mbox{Li}_2(1/u_0)+\alpha(z_0,u_0))
=2.82812208...
\eeq
in agreement with \cite{AHW}.

\subsection{Knot $6_1$}

The sum (\ref{6_1}) at large $N$ is approximated by the triple integral
\begin{eqnarray}\label{contour6_1}
\langle6_1\rangle&\sim&\int dzdudv\  
\exp\frac{i}{2\gamma}[
\mbox{Li}_2(z)-\mbox{Li}_2(1/z)-\mbox{Li}_2(u)
+\mbox{Li}_2(1/v)\nonumber\\
&+&\alpha(uv/z,z/u)+2\pi i\log(u/z)],
\end{eqnarray}
where $\alpha(z,u)$ is defined in (\ref{alpha}). The stationary points are
solutions to the algebraic system of equations:
\beq
z(1-z)^2=-u^2v,\quad z^2(1-u)=u^2v,\quad z(1-v)=-uv.
\eeq
The maximal contribution to the integral comes from the solution 
$(z_0,u_0,v_0)$ with
\beq
\mbox{Im}\ z_0<0<\mbox{Im}\ (u_0v_0).
\eeq
Thus, the asymptotics of the absolute value of the integral reads
\beq
|\langle 6_1\rangle|\sim \exp\frac{V(6_1)}{2\gamma},\quad \gamma=\pi/N\to0,
\eeq
where
\begin{eqnarray}
V(6_1)=-\mbox{Im}\ (\mbox{Li}_2(z_0)-\mbox{Li}_2(1/z_0)-\mbox{Li}_2(u_0)
+\mbox{Li}_2(1/v_0)\nonumber\\
+\alpha(u_0v_0/z_0,z_0/u_0)+2\pi i\log(u_0/z_0))=3.16396322...
\end{eqnarray}
again in agreement with \cite{AHW}.

\section{Summary}

We have demonstrated on three examples of hyperbolic knots, that the
link invariant, defined in \cite{K1,K2} via cyclic quantum dilogarithm,
in the asymptotic limit $N\to\infty$ reproduces the hyperbolic volume
of a knot, see formula (\ref{main}). This result implies a possible relation
of the corresponding combinatorial TQFT to quantum 2+1-dimensional gravity.

\section{Acknowlegements}

The author is grateful to M. Blau, L. Freidel, J.M. Maillet, F. Smirnov for 
valuable discussions, and L.D. Faddeev for his encouragement in this work. 
The work is supported by the Programme TEMPRA-Europe de 
l'Est from the R\' egion Rh\^ one-Alpes.

\end{document}